# Investigations of complex systems' dynamics, based on reduced amount of information – part 3 – Dynamical Phenomena Indicator – the fastest and most universal approach for investigations of dynamics of complex networks of coupled nonlinear systems

Volodymyr Denysenko, Artur Dabrowski

Division of Dynamics, Lodz University of Technology,
ul. Stefanowskiego 1/15 Lodz, Poland.
Volodymyr.Denysenko@dokt.p.lodz.pl
Ar2rDe@p.lodz.pl

**Abstract**

Recently, we have demonstrated that our approach is a highly effective tool while analysing complex phenomena existing in networks of coupled nonlinear systems. In the present article we present the results of our investigations into a specific aspect of the presented method. We prove its effectiveness while applying for fast investigations of complex systems and easy detection of different uncommon dynamical phenomena states. We also extend our method introducing new Dynamical Phenomena Indicator (DPI), designed especially for effective detection of complex dynamical phenomena states in the wide range of the parameters of complex networks of coupled nonlinear systems.

Contrary to commonly applied methods, the proposed approach allows for identification of complex dynamical phenomena long before stabilization of the system. The method bases on early signalized tendency of the system to split its dynamics to separately synchronized subsystems.

The most important fact is that proposed approach is highly universal and can be applied for both, symmetrical and non-symmetrical topologies of coupling as well as networks of identical and non-identical oscillators.

Moreover, since DPI values are obtained from the current state of dynamical system given by values of the system variables, proposed method of fast searching has a huge potential for experimental application. Following this reasoning the presented results can be treated both as exemplary numerical investigations and analysis of experimentally obtained results.

## 1. Introduction

Collective behaviour of complex networks(CN) of coupled oscillators attracts attention of many scientists due to existence of many complex patterns of their dynamic behaviour, as well as applications in many different disciplines of science and engineering. Even that the problem was investigated for many decades, there is still lack of universal, efficient tools allowing for recognizing patterns existing in those, the most complicated systems. Complexity of topologies of CN influences studies of such systems to be the most complicated and time demanding cases. Simultaneously what attracts scientist for decades in dynamics of complex networks(CN) of coupled oscillators is mainly huge spectrum of different types of collective behaviours.

The first of them, probably the most significant is synchronization, which was for the first time observed in 1673 by Huygens [1] for the case of synchronization of two clock pendulums. Since then, knowledge about this state has been developed by many scientists and applied in various engineering fields, such as: electrical systems [2], mechanical [3], [4], [5], their combination, electro-mechanical systems [6], cryptography [7] and other branches of science [8]. It was also was used to explain energy flows in dynamical systems [9] or different nature phenomena, such as extinction or survival of species [11], or synchronization in complex social human networks [12]. It was observed in processes in human brain [10], where it was found exemplarily that synchronization of neurons may cause pathologies in brain function, such as Parkinson disease [13] or epilepsy [14]. Understanding of this process can help us to predict seizures [15] or early diagnose Parkinson's disease [16]. On the other hand dynamics of such a complex systems like brain is so complex that the same synchronization phenomena constitutes basis of learning processes [17].

Another unexplained and unexpected collective CN behaviour was later discovered in one of the simplest types of coupled networks: a network consisting of symmetrically coupled identical oscillators. Due to the symmetry of the considered system, it seems obvious that all the oscillators can either synchronize in some symmetrical patterns or just drift incoherently. However Y. Kuramoto [18] discovered that under some specified conditions this symmetry breaks and coherent and incoherent oscillators can coexist in such networks. This unexpected phenomenon was called "chimera state" by Steven Strogatz in [19]. Since then, such behaviour has been found in many different models, such as: mechanical and electrical systems [20] [21], chemical oscillators [22], memristive neuron networks [23] and others. Various patterns in the dynamics of complex networks (CN) are still being discovered, such as traveling chimera state, where the synchronized segment can migrate in a stable form, or the so-called state or pulse (breathing chimera).

Despite significant progress, numerous challenges to Complex Networks analysis persist even today. Many of them concern developing new methods of detecting interesting behaviours, and designing approaches on CN control [25] [26].

A few different indicators and methods have been proposed for detecting complex dynamical phenomena (CDP). Most commonly used indicator, referred to as Strength of Incoherence (SI), was introduced in [27]. This indicator involves dividing oscillators into groups and calculating the standard deviation of the differences in their dynamical states, averaged over time for each group. While this approach provides information about the presence of a chimera state given certain initial conditions and system parameters, it does not reveal the behaviour pattern of the oscillators or their tendency to synchronize. Another indicator, called the local order parameter, which can be also utilized for detecting of grouping oscillators in different patterns, was presented in [28]. Additionally, calculating the average velocities of oscillators [29] or their FFT [30] [31] can also indicate the existence of these states and provide insights into oscillators dynamics. However, all these methods share a common drawback: they only give useful information about the investigated CN after it has stabilized. This limitation is absent in our presented method. This approach bases on recently presented, novel method for studying dynamics of CN, namely Transversely Directed Lyapunov Exponents (TDLE) [32].

It is worth to note that TDLE are obtained from the current state of investigated dynamical system. As a result, this method can be successfully applied for both numerical simulations and experiments. It also influences on universality of the approach which can be applied in any type of coupled systems, regardless of type of coupling and differences or similarities between oscillators.

However, the most significant feature of the TDLE approach is that it provides information about the early tendency of CN dynamics long before stabilization occurs. This feature was leveraged in order to develop a new, efficient algorithm for CN investigation. Presented algorithm accelerates scanning of the system parameters space and initial conditions for potentially interesting dynamic behaviours, making it arguably the most effective method for CN research.

## 2. Introduction to DPI approach

### 2.1 TDLE method

The main purpose of the TDLE method is to analyse system behaviour in directions transversal to the synchronization manifold.

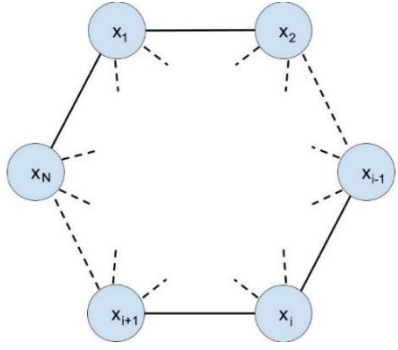

*Fig. 1 Symmetrically coupled identical oscillators*

Assume that a CN consists of *n* symmetrically coupled identical oscillators (Fig. 1). The system can be described by the following set of coupled differential equations:

$$\dot{x} = f(x,t) + \alpha(G \otimes H)x \qquad (1)$$

where α is the coupling coefficient, **G** is the connectivity matrix, **H** is the matrix that chooses variables engaged in the coupling, $\otimes$ is the direct Kronecker product of two matrices. Eq. (1) can be rewritten in form:

$$\frac{d x_j}{dt} = f(x_j, t) + g_j(x_1, x_2, \ldots, x_n), j = 1 \ldots n \qquad (2)$$

where $f(x_j, t)$ are the functions describing each of the *n* nodes, $g_j(x_1, x_2, \ldots, x_n)$ are the functions describing couplings of these nodes, that depends on values of α, **G** and **H**, and $(x_1, x_2, \ldots, x_n \in R^m)$ are vectors of state variables of each of nodes.

The differences among combinations of the state variables can be introduced as upper triangular matrix of perturbation vectors $z_{ij}$ introduced in directions transversal to synchronization manifold:

$$Z = \begin{bmatrix} 0 & z_{12} & \cdots & z_{1n} \\ 0 & 0 & \cdots & z_{2n} \\ \vdots & \vdots & \ddots & \vdots \\ 0 & 0 & \cdots & 0 \end{bmatrix} = \begin{bmatrix} 0 & |x_1 - x_2| & \cdots & |x_1 - x_n| \\ 0 & 0 & \cdots & |x_2 - x_n| \\ \vdots & \vdots & \ddots & \vdots \\ 0 & 0 & \cdots & 0 \end{bmatrix} \qquad (3)$$

Presented TDLE approach bases on the idea of commonly applied Lyapunov Exponents (LE). However contrary to LE approach, where directions of the perturbation vectors evolve in time and are orthogonal to the trajectory of the system, perturbations introduced TDLE method have pre-selected directions transversal to the synchronization manifold (SM). Each of the vectors $z_{ij}$ span *m*-dimensional perturbation subspace $S_{i,j}$, transversal to the SM, where *m* is dimension of each of the node system. In each of these subspaces $S_{i,j}$, an evolution of perturbation can be analysed separately using exponential averaging, similar to Lyapunov exponents and obtaining thereby the values of Transversely Directed Lyapunov Exponents (TDLEs) see Eq.(4). Note that thanks to this construction, no orthogonalization has to be introduced to the algorithm, unlike in the case of the commonly applied Lyapunov exponents (LE).

$$Z(t) = \begin{bmatrix} 0 & |z_{12}(t)| & \cdots & |z_{1n}(t)| \\ 0 & 0 & \cdots & |z_{2n}(t)| \\ \vdots & \vdots & \ddots & \vdots \\ 0 & 0 & \cdots & 0 \end{bmatrix} = \begin{bmatrix} 0 & |z_{12}^0|e^{DLE_{12} \cdot t} & \cdots & |z_{1n}^0|e^{DLE_{1n} \cdot t} \\ 0 & 0 & \cdots & |z_{2n}^0|e^{DLE_{2n} \cdot t} \\ \vdots & \vdots & \ddots & \vdots \\ 0 & 0 & \cdots & 0 \end{bmatrix} \qquad (4)$$

where $|z_{i,j}^0|$ are respective initial perturbations in each of the introduced subspaces transversal to the synchronization manifold.

The values of TDLE indicators can be obtained using any method of Lyapunov exponents estimation, for example one presented in [35]. This article makes use of combination of two methods designed by authors of the article, the method of extracting LE from differences between dynamical states of two coupled systems [33] and fastest and simplest method of LE estimation for continuous time dynamical systems analysis as presented in [34].

### 2.2 Application of the special behaviour of TDLEs in complex dynamical phenomena (CDP) analsis – spreading of TDLEs

According to the introduced definition of TDLEs, the negative values of TDLEs mean that corresponding pairs of oscillators tend to synchronize, while the positive values of TDLE mean that corresponding pairs of oscillators tend to diverge. However, due to limitation of attractor volume, TDLEs for incoherent pair of oscillators have values that are close to zero, not positive as expected. Since it does not have a significant qualitative impact on the presented results and conclusions, corrections to the algorithm will be developed in the future investigations.

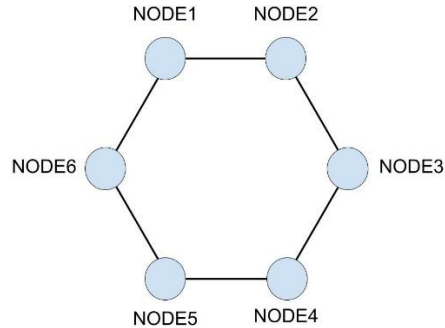

Fig. 2 Schematic topology of 6 periodically forced Duffing locally coupled oscillators

To give a better understanding of TDLE behavior, let's consider the following example: a CN that consists of 6 periodically forced Duffing nonlocally coupled oscillators.

$$\begin{pmatrix}\dot{x}_{2j}\\ \dot{x}_{2j+1}\end{pmatrix} = \begin{pmatrix}x_{2j+1}\\ -2hx_{2j+1} - kx_{2j}^3 + F\cos(\omega t)\end{pmatrix} + \sigma \begin{bmatrix}-2 & 1 & 0 & 0 & 0 & 1\\ 1 & -2 & 1 & 0 & 0 & 0\\ 0 & 1 & -2 & 1 & 0 & 0\\ 0 & 0 & 1 & -2 & 1 & 0\\ 0 & 0 & 0 & 1 & -2 & 1\\ 1 & 0 & 0 & 0 & 1 & -2\end{bmatrix} \otimes \begin{bmatrix}0 & 0\\ 1 & 0\end{bmatrix}\begin{pmatrix}x_{2j}\\ x_{2j+1}\end{pmatrix} \quad (5)$$

j = 0..5

Fig.(2) shows such coupling schematically, where each node represents one Duffing oscillator, and each line represents linear coupling between oscillators. Figs.(3-5) show time series of TDLE spectrum for different dynamical states of the network.

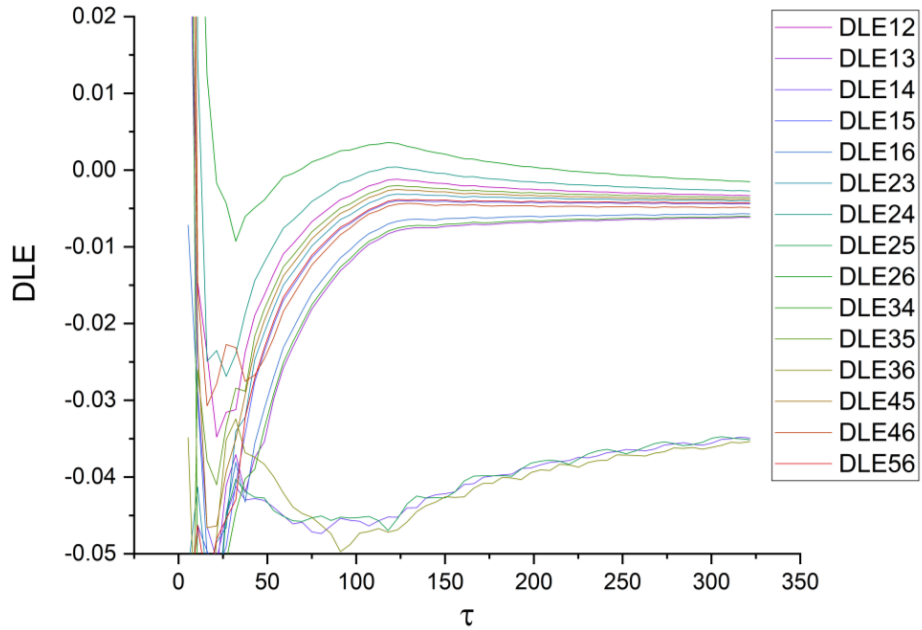

*Fig. 3 Fast separation of TDLE values for 6 locally coupled Duffing oscillators*

Fig.(3) shows the case, where synchronization of nodes 1-4, 2-5 and 3-6 can be detected. It can be noticed that after approximately 325 periods of oscillations TDLE values stabilize. It is worth noting that a significant difference between negative TDLE values and TDLE values that are close to zero appears after just 50 periods of oscillations.

Similar phenomena can be observed in Fig.(4) – however TDLEs spread in this case later, and before 300 periods of oscillations and the differences between values of spectrum is almost unnoticeable.

Fig.(5) depicts case of incoherent oscillations. It is worth noting that the difference between values of TDLE spectrum in such case is unnoticeable. So, we can conclude that for any possible CDP, there appear differences between values of TDLEs, corresponding to different pairs oscillators, based on its tendency to synchronize. These differences can be noticed before stabilization of the system and synchronization of nodes. In the article it is proposed to utilize this phenomenon to obtain algorithm for fast scanning of parameters space, while searching of possibility of complex dynamical behaviours. While comparing Fig.(3) and Fig.(4) it can be seen that for network with the same oscillators and topology achieving of TDLEs stabilization and TDLEs spreading differs much. As the main task of our investigations was to develop a fast, time-efficient algorithm, this fact was the main challenge in deriving universal and optimized parameters of the algorithm.

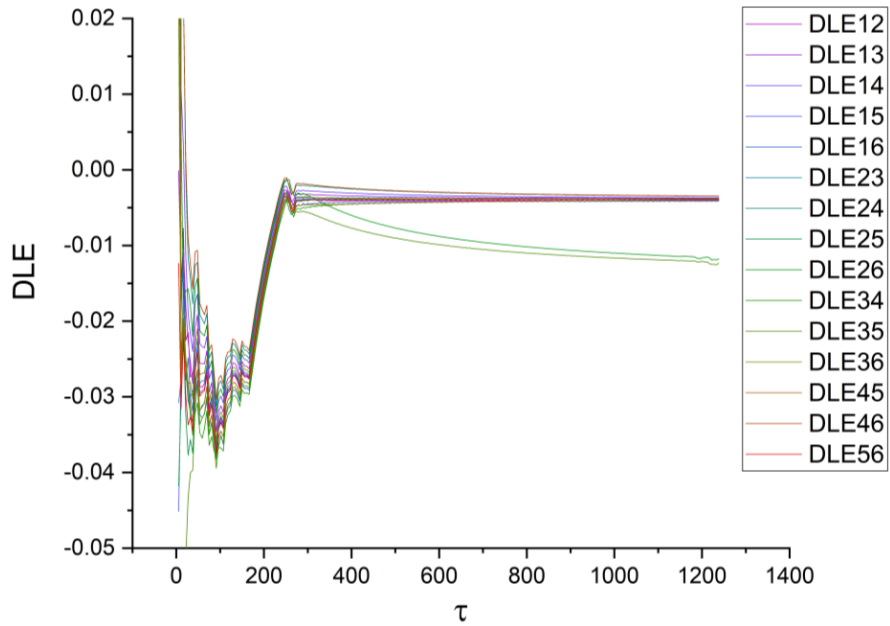

*Fig. 4 Slow separation of TDLE values for 6 locally coupled Duffing oscillators*

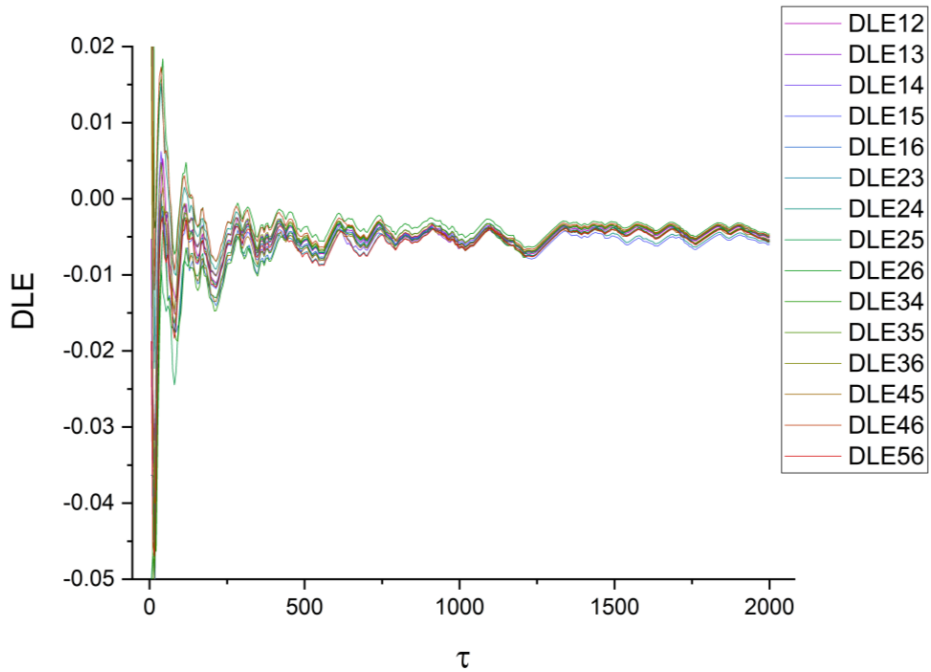

*Fig. 5 TDLE spectrum for incoherent oscillations for 6 locally coupled Duffing oscillators*

**2.3 Proposed algorithm of complex dynamical phenomena recognition**

Presented feature of early spreading TDLEs in case of complex dynamical phenomena existence is proposed to be utilised to develop the most effective and universal way of CDP detection. Analysing other options, such as the spread of norms of differences between dynamical states, was excluded, as it was considered to be less stable way of detecting the aforementioned spread. For the concerned case, where dynamics of the system partly approaches synchronization, norms of differences between dynamical states is exponentially lowering to zero. It is connected with exponentially variable mean value and more complicated to include in the algorithm condition of stabilization, unlike to TDLE measure striving for constant value. The spread of average frequencies applied in case of chimera states existence was also excluded as not universally efficient condition of CDP detection in complex networks(CN), especially when analysing dynamics divided into parts with not so evident differences of frequencies, what is the often appearing case in such a systems.

To utilize the previously described feature of TDLE, the following algorithm was proposed: TDLE values were calculated and ordered from largest to smallest. Then, the differences between consecutive values from the larger TDLE value to the smaller TDLE value were calculated. The maximum obtained difference ($\Delta DLE_{max}$) was compared to threshold ($\Delta DLE_{thres}$) which has to be experimentally determined. If this maximum difference exceeds the threshold, it is interpreted as a possible indication of a CDP for that parameter set. If the maximum difference remains below the threshold for a sufficiently long time, it can be interpreted as a lack of the possibility of a complex state.

As it was mentioned, while analysing examples depicted in Fig.(3) and Fig.(4), there can be observed different scenarios of synchronization of particular couples of oscillators. In the case, shown in Fig.(3), gap between negative and nonpositive values of TDLE is easily noticeable and appears after only 50 periods, what can be interpreted as strong tendency to synchronize. For the case, depicted in Fig.(4), the difference between negative and nonpositive values of TDLE appears almost after 400 periods and then slowly increases in time. To consider such possibilities and obtain universal algorithm, threshold value must be not constant value, but a function of time.

To find appropriate threshold function $\Delta DLE_{thres}$ with minimal effort, it is proposed to use the data from the simulations for CN system with less number of oscillators than investigated CN. Fig. 6 shows time series of maximum differences $\Delta DLE_{max}$ for different coupling coefficients and all symmetrical topologies of 4 Duffing oscillators in cases of detected incoherent oscillations. Additionally, three obtained threshold functions are introduced in the Fig.(6). Values of the $\Delta DLE_{thres}$ in each moment of time were obtained in such a way as to have a value greater than defined percentage of obtained $\Delta DLE_{max}$ values for incoherent oscillations. To differentiate possible complex dynamical states from uncoherent states, the threshold function must be approximated to be close to the largest value of $\Delta DLE_{max}$ obtained from numerical simulations at any given time. Then resulting threshold function can be applied to larger networks of the same oscillators.

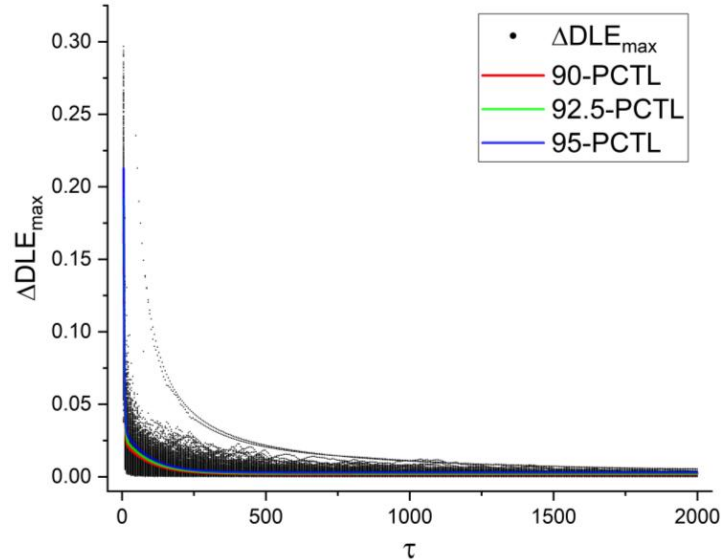

*Fig. 6 Time series of maximum TDLE differences for various coupling coefficients and all symmetrical topologies of 4 Duffing oscillators.*

In the Fig.(7) there are presented obtained functions $\Delta DLE_{thres}$ together with value of $\Delta DLE_{max}$ presented in Fig.(3). It can be observed that value of $\Delta DLE_{max}$ crosses the $\Delta DLE_{thres}$ function after approximately 60 periods of oscillations; therefore, it can be concluded that a complex dynamical state appears in this case.

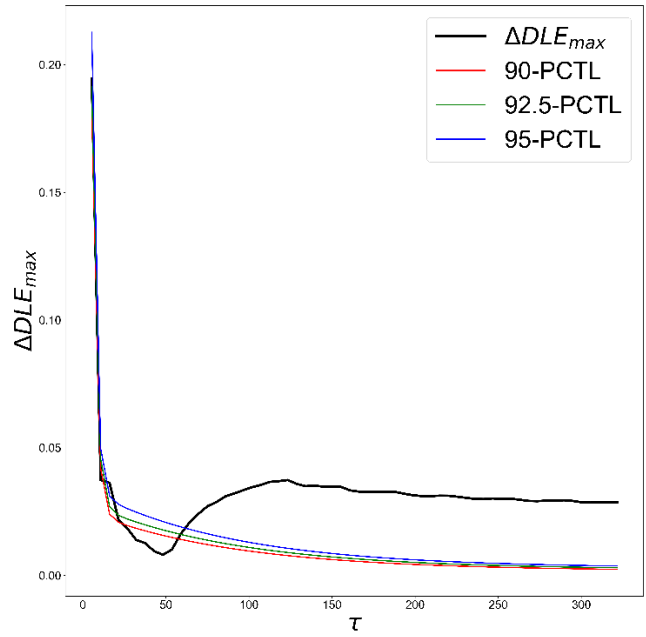

*Fig. 7 Fast TDLE separation detection example for 6 locally coupled Duffing oscillators*

Similarly in Fig.(8), the obtained $\Delta DLE_{thres}$ functions are presented alongside the $\Delta DLE_{max}$ values shown earlier in Fig.(4). It can be observed that $\Delta DLE_{max}$ crosses the $\Delta DLE_{thres}$ function after approximately 500 oscillation periods, indicating the appearance of a CDP in this case.

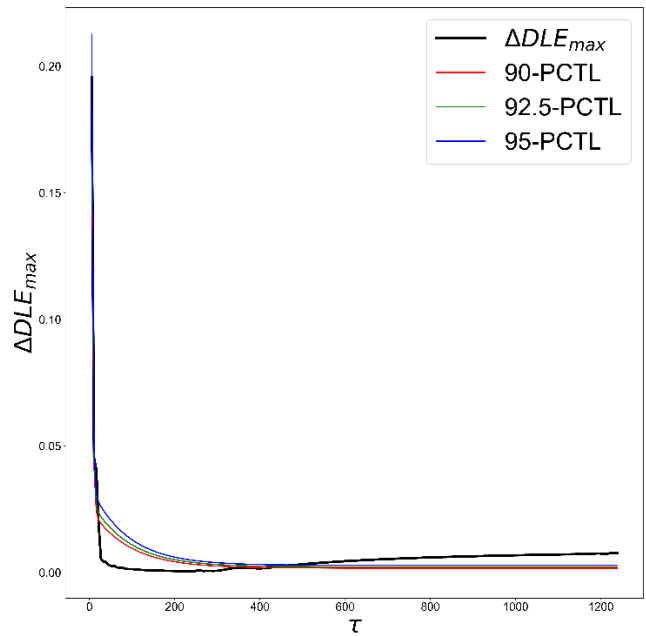

*Fig. 8 Slow TDLE separation detection example for 6 locally coupled Duffing oscillators*

## 2.4 Master Stability Function

As application of Transversal Lyapunov Exponents (TLE) and the Master Stability Function (MSF) [36] is the most effective tool for identification of complete synchronization states of sets of coupled identical oscillators it is proposed to apply this tool in the beginning of investigations. Then excluding these ranges of complete synchronization allows to easily find regions where it is possible for the system to exhibit more complex dynamical patterns.

To introduce MSF approach the block form of differential equations (Eq.(1)) has to be used. Then tendency of a network of oscillators to synchronize/desynchronize is quantified by the eigenvalue spectrum $\lambda_i$ ($i=1..n/2$) of the connectivity matrix G, see Eq.(5), i.e. the Laplacian matrix representing the topology of connections between the network nodes. After block diagonalization of the variational equation of Eq.(1) there appear $i$ separated blocks:

$$\dot{\mathbf{z}}_i = (D\mathbf{f} + \alpha\lambda_i D\mathbf{H})\mathbf{z}_i \tag{6}$$

where $\mathbf{z}_i$ represent different modes of the perturbation with respect to the synchronization manifold. Eq.(6) allows for calculation of the largest Transversal Lyapunov Exponent (TLE) presented in Fig.(11) over variable coupling coefficient $\alpha$.

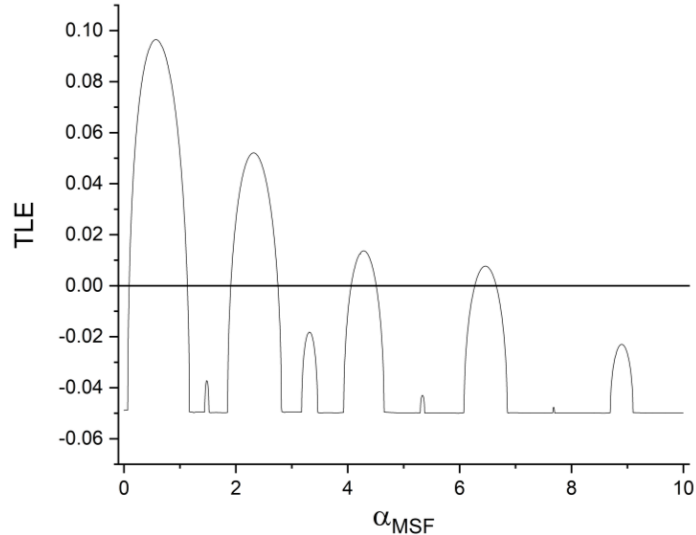

*Fig. 9 Transversal Lyapunov exponent.*

Having regard to Eq.(6), chart in Fig.(11) has to be rescaled in the horizontal - $\alpha$ direction, with rescaling eigenvalues $\lambda_i$. After this procedure one obtains number of i-charts which together stand as Master Stability Function (MSF). The example of MSF for CN, described in previous chapter (Fig. 2) is shown in Fig.(12). The MSF providing all negative TLEs defines stable ranges of synchronization manifold, while for remaining ranges other types of CN dynamics occur. Approach proposed in the article assumes conducting investigations in the ranges where there were detected at least one positive TLE.

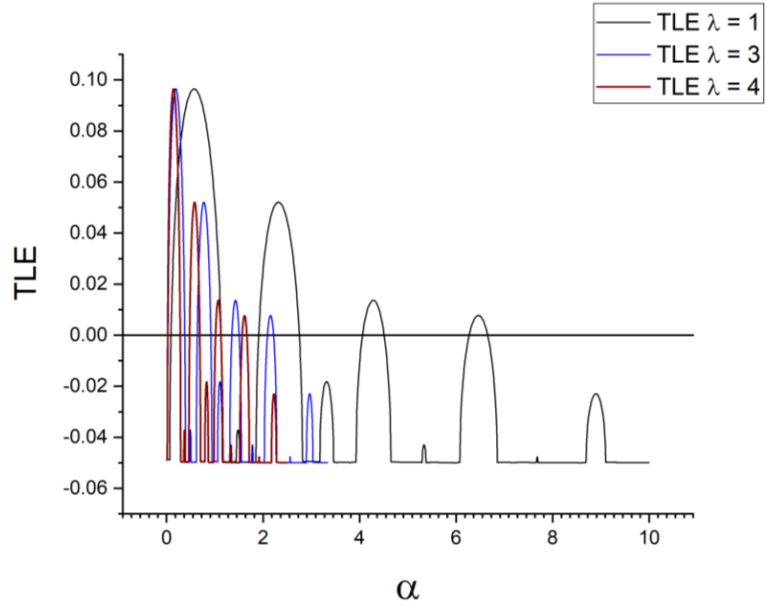

*Fig. 10 MSF for locally coupled ring of 6 Duffing oscillators(Fig. (2))*

## 3. Methodology and numerical investigations

The first step of experiments involved primary identification and excluding ranges of coupling strength coefficient where complete synchronization occurs, using Transversal Lyapunov Exponents (TLE) and the Master Stability Function (MSF) [36]. Then, initial numerical simulations were carried out on the ranges of parameters that were excluded in the first step, where the considered system potentially exhibits interesting dynamical behaviour. It was carried out on networks that consisted of 4 to 8 oscillators with all possible types of symmetrical couplings. All identified ranges of coupling strength and random initial conditions were chosen for test, to assure presence of different types of behaviour of CN, including fast synchronization of particular nodes, slow synchronization of particular nodes and incoherent states. Initial simulations included perturbation vector analysis and calculation of TDLE's. Time series of perturbation vectors and TDLE's obtained during initial simulations were saved and then their analysis was applied for threshold function approximation, for further algorithm evaluation and its validation.

To approximate threshold function, required for proposed algorithm, data obtained from numerical simulations for 4 oscillators with both local and global symmetrical coupling was analysed as described in paragraph 2.3.

Then, algorithm was validated analysing data for systems, consisting of 4 to 8 oscillators with all possible types of symmetrical couplings.

In initial simulation algorithms, the selection of termination conditions is crucial when detecting a CDP. To streamline the procedure and to avoid introducing unnecessary calculations to the algorithm, it appears reasonable to perform CN simulation until one of following conditions will be satisfied: synchronization between at least single pair of oscillators will occur, or simulation will exceed the predefined time limit(set as 2000 periods of oscillations). After initial studies the condition of synchronization was defined as the case when value of length of perturbation vector **z** is less than $10^{-15}$ as the easiest to implement. If condition was satisfied in the defined time range, corresponding set of initial conditions and parameter value was denoted as one, for which probability of complex state occurrence was high. Data obtained in that way was used for proposed algorithm tests.

In order to measure the stabilization of the TDLE value, a buffer of a selected fixed size was defined and introduced. In this investigations, a buffer capacity of 200 elements was chosen as a practical compromise between TDLE estimation accuracy and simulation duration. Following each computational step, the current TDLE value was stored in the buffer. Once the buffer reached its maximum capacity, the standard deviation of all TDLE values within the buffer was computed. If the standard deviation fell below a predetermined threshold, indicating stable TDLE values, the calculations were concluded. If not, the buffer was cleared, and the process was repeated. This approach was applied for TDLE's calculating during whole simulations.

All the programs for conducting numerical simulations have been written in C++ by means of the Visual Studio Code environments. The Runge-Kutta method of the fourth order (RK4) has been used to solve ordinary differential equations. The integration step has been adjusted for each analysed system individually, taking into account its own time scale. Scripts for analysis and validation of proposed algorithm were written in Python. Origin software was used for curve fitting, what is necessary for threshold function approximation.

**Coupling radius**

In order to differentiate topologies of investigated CN, coupling radius R was introduced. Meaning of coupling radius is as follows: it is a number of couplings with neighbouring oscillators divided by 2 as CN is symmetrical. Coupling radius must be taken into account during connectivity matrix construction, as described by Eq.(7).

$$G_{ij} = \begin{cases} \max(-(N-1), -2R), if\ i = j \\ 1, if\ |i-j|\ mod\ N \leq R \\ 0, if\ |i-j|\ mod > R \end{cases} \quad (7)$$

where N is a number of oscillators in system and R is a coupling radius.

**3.1 Introduction of the new Dynamical Phenomena Indicator (DPI) for determination of the dynamical states of complex networks of coupled nonlinear oscillators**

To describe and differentiate dynamical state patterns that occur in system, based on results of TDLE calculations, the Dynamical Phenomena Indicator (DPI) was introduced:

$$DPI = \frac{N_{aver}}{N_{unsync} + N_{pairs}} \quad (8)$$

where:

- ✓ $N_{aver} = \frac{N_{sync}}{N_{groups}}$ – average number of synchronized pairs per group of synchronized oscillators. The number provides growing indicator values together with growth of the dynamical structures of synchronized oscillators, meaning simultaneously growth of complexity of dynamics of the system,
- ✓ $N_{sync}$ – number of pairs of synchronized oscillators,
- ✓ $N_{Groups}$ – number of separate groups of synchronized oscillators,
- ✓ $N_{pairs}$ – maximal number of synchronized pairs of oscillators for analysed network consisted of $N$ coupled oscillators. The number provides DPI normalization growing from zero in case of all oscillators coherent to value 1 for the case of complete synchronization,
- ✓ $N_{unsync}$ – number of pairs of unsynchronized oscillators provides DPI uniqueness. The bigger is the number the less complex is the system configuration.

Such indicator allows us to recognize dynamical state of the system, differentiate coherent and uncoherent states and to make conclusions on complex patterns that occur in a system. Examples of different patterns together with DPI values for 8 coupled oscillators are shown in Fig.(11), where oscillators are represented with blue dots, and synchronised oscillators are connected with green lines. Note that examples are organized in a way that DPI values grow up from the left to the right side together with sizes of synchronized structures.

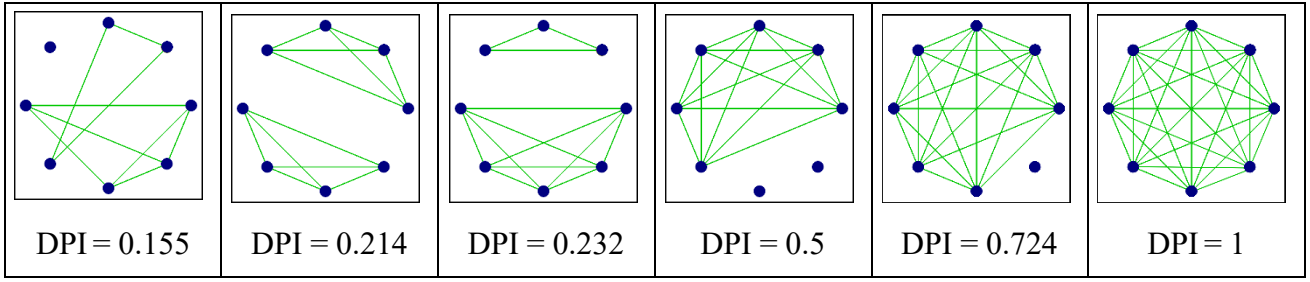

| DPI = 0.155 | DPI = 0.214 | DPI = 0.232 | DPI = 0.5 | DPI = 0.724 | DPI = 1 |

*Fig. 11 Examples for 8 coupled oscillators with different configurations of synchronized nodes*

Dynamical behaviour of the networks, observed during investigations, can be divided into two general cases. The first one is coexistence of groups consisted of the same number of synchronized oscillators, and the second one is coexistence of groups of two different numbers of synchronized oscillators. Those two cases are familiarised beneath, where we present theoretical analyses of DPI index values and compare them with results of configurations found during conducted simulations. We also analyse DPI uniqueness to prove that it gives both information about dynamical state of the system and possibility to recover its dynamical configuration.

Let us firstly consider the case, where dynamics a of $N$ coupled oscillators is divided into $N_{Groups}$ separate groups, each consisted of equal number $N_G$ synchronized oscillators. Then, from the basics of combinatorial analysis the maximal number of synchronized pairs in one group:

$$N_{pairs/group} = \frac{N_G!}{2!(N_G-2)!} = \frac{(N_G-1)N_G}{2} \quad (9)$$

Then number of all the synchronized pairs:

$$N_{Sync} = N_{Groups} N_{pairs/group} \quad (10)$$

and:

$$N_{unsync} = N - N_{Sync} \quad (11)$$

Simultaneously, maximal number of synchronized pairs of oscillators for considered system consisted of N coupled oscillators:

$$N_{pairs} = \frac{N!}{2!(N-2)!} = \frac{(N-1)N}{2} \quad (12)$$

Finally substituting to Eq.(8):

$$DPI = \frac{(N_G-1)N_G}{2(N - \frac{2N_{sync}}{N_G-1} + \frac{(N-1)N}{2})} \quad (13)$$

Dependence of DPI values for $N=100$ oscillators is presented in Fig.(12) as the function of $N_{sync}$ - number of pairs of synchronized oscillators and $N_G$ – number of synchronized oscillators per group. It can be seen that presented dependence is monotonically unique. It means that configuration of the synchronized structures can be recovered from DPI values.

It can be also seen that DPI values grow up together with growing number of oscillators in one group, achieving value equal to 1 for the case of complete synchronization. Note that presented dependence is continuous opposite to discrete values of DPI for real system consisted of discrete numbers of oscillators which are presented in Fig (13).

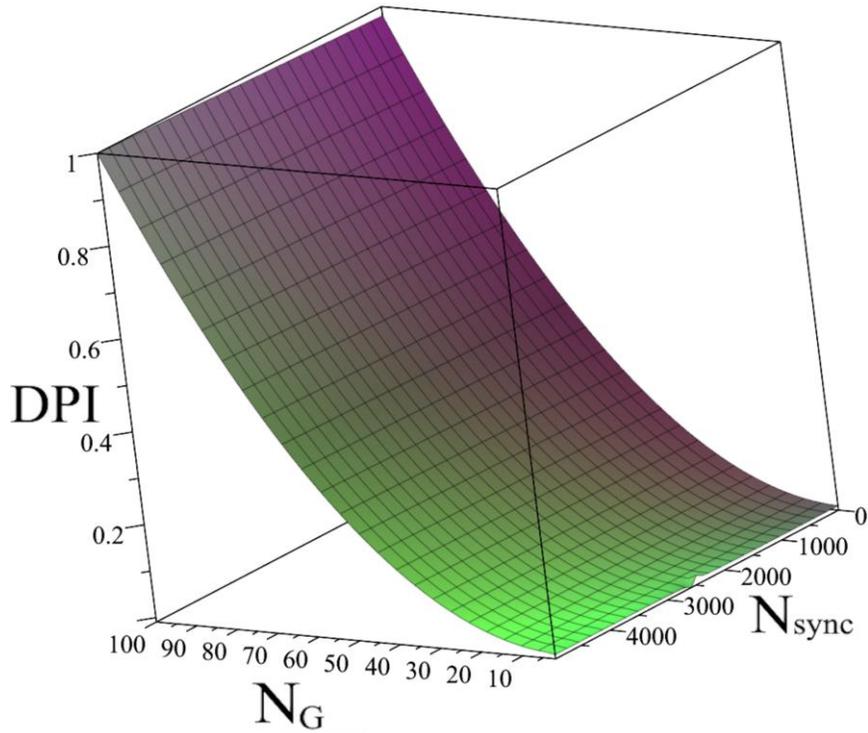

*Fig. 12 Dependence of DPI on $N_{sync}$ - number of pairs of synchronized oscillators and $N_G$ – number of synchronized oscillators per group*

Fig (13) presents comparison of the analytically obtained continuous values of DPI (continues curves) and discrete values (star-marked), obtained for configurations found during numerical investigations (Table 1). CN system that consisted of 8 coupled oscillators was taken under consideration. In the Fig (13) analytical functions (cross sections of Fig.(12) are presented for groups consisted of equal numbers of $N_G$ oscillators changing from 2 to 8. In this figure there can be seen agreement between theoretical and simulation results.

There can be also noted that results values of DPI grow monotonically together with growing number of synchronized oscillators in each group and number of synchronized pairs.

While analysing numbers of star-marked points for each case it can be seen that they fulfil all the possible cases of configurations in the case of 8 coupled oscillators. For example while groups consist of $N_G$ equal to 2 oscillators there are 4 possible cases – 1or 2 or 3 or 4 existing groups of synchronized oscillators. Accordingly in case of $N_G$ equal to 3 and 4 oscillators there are 2 possible cases. Finally from the $N_G$ equal to 5 there is possible only one case. It suggests that in the case of 8 coupled Duffing oscillators no special dynamics pushing system to some special dynamical direction occurs.

Finally note that each of the continuous curves finishes with star-marked configurations detected from numerical simulations. These points concern cases of $N_G$ equal to subsequent values of maximal numbers of synchronized pairs in CN consisted of corresponding numbers of oscillators. Thus value $N_{sync}$ =28 corresponds to DPI= 1 and depicts full synchronization state.

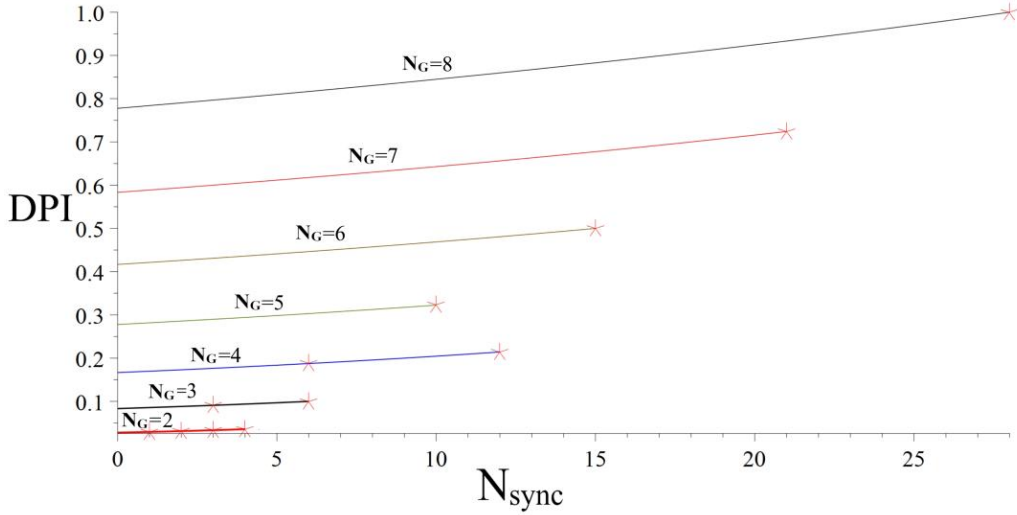

*Fig. 13 Dependence of DPI on $N_{sync}$ - number of pairs of synchronized oscillators and $N_G$ – number of synchronized oscillators per group*

Let us now present analysis of the second observed case, when it was noticed coexistence of two types of groups consisted of two different numbers of synchronized oscillators. For an assumed CN topology, dynamical configuration state of this system with two types of synchronized groups is given by following parameters: numbers of separately synchronized groups $N_{Groups1}, N_{Groups2}$ of the first and second type, and numbers of oscillators engaged in the first and second type of groups $N_{G1}, N_{G2}$. Then number of not synchronized oscillators:

$$N_{unsync} = N - (N_{Groups1} \cdot N_{G1} + N_{Groups2} \cdot N_{G2}) \qquad (14)$$

It can be seen from Eq(16), that DPI value depends on five parameters. Full analysis of all the possible cases would be very extensive and would not bringing any special benefits. Especially that main idea of presented in the article approach is to find a way of fast scanning the space of system parameters for the purpose of indicative determination of possibility of some special dynamical behaviors existence. Then the interesting areas that have been discovered can be explored using any other methods.

Thus to simplify analysis and presentation of the DPI dependence on growing up synchronized structures, we assumed simple dependences between numbers of groups $N_{Groups1}, N_{Groups2}$ and also between numbers of oscillators engaged in the first and second type of groups $N_{G1}, N_{G2}$.

Firstly we assumed that these two types of groups of synchronized oscillators differ in size by one oscillator and for the independent variable numbers of oscillators engaged in the first type of groups $N_{Groups1}$ was taken. Then, dependent on $N_{Groups1}$, the number of oscillators engaged in the second type of groups equals to:

$$N_{Groups2} = N_{Groups1} + 1 \qquad (15)$$

For further simplification we assumed that the system consists of constant number $N=100$ of coupled oscillators. Additionally we also assumed that number of existing both types of groups was chosen in a way to engage maximally possible number of oscillators. Thus together with growth of number of first groups, the number of the second types was lowering engaging other oscillators not involved in the first structures. In that way we obtained one parameter $N_{Groups1}$ identifying size of appeared structures.

Assuming given above limitations and introducing them to Eq.(8) the final formula for DPI in this case has the following form:

$$DPI = \frac{N_{Groups1}(N_{G1}-1)N_G + \left\lfloor \frac{-N_{Groups1} N_{G1}+100}{N_{G1}+1} \right\rfloor N_{G1}(N_{G1}+1)}{2\left(N_{Groups1}+\left\lfloor \frac{-N_{Groups1} N_{G1}+100}{N_{G1}+1} \right\rfloor\right)\left(-N_{Groups1}N_{G1}-\left\lfloor \frac{-N_{Groups1} N_{G1}+100}{N_{G1}+1} \right\rfloor(N_{G1}+1)+5050\right)} \qquad (15)$$

where:
- $N_{G1}$ - number of synchronized oscillators in each first-type group,
- $N_{Groups1}$ - number of first type of groups consisted of $N_{G1}$ synchronized oscillators,
- $\lfloor x \rfloor$ – integer part of x.

Note that similarly to earlier analysis of one type of groups coexistence, presented dependence is continuous showing general tendency of changes of DPI values together with change of continuously growing up sizes of synchronized structures, which in real constitute discrete set of points immersed in an infinite set of points constituting a surface presented in Fig.(14). It the Figure there is visible growth of DPI values together with growing number of synchronized oscillators in first groups, meaning also growing number of synchronized oscillators in the second groups, following Eq.(14). It can be also seen from the red color line connecting black points that together with growth of number of groups, meaning groups are smaller, DPI also is lowering down. Note that the number of black marked picks of DPI values is equal to number of transitions between numbers of groups of synchronized oscillators from one to ten groups. Monotonically lowering down values of DPI for black marked points are expected values for real discrete DPI values, proving uniqueness of the presented indicator for assumed number of coupled oscillators.

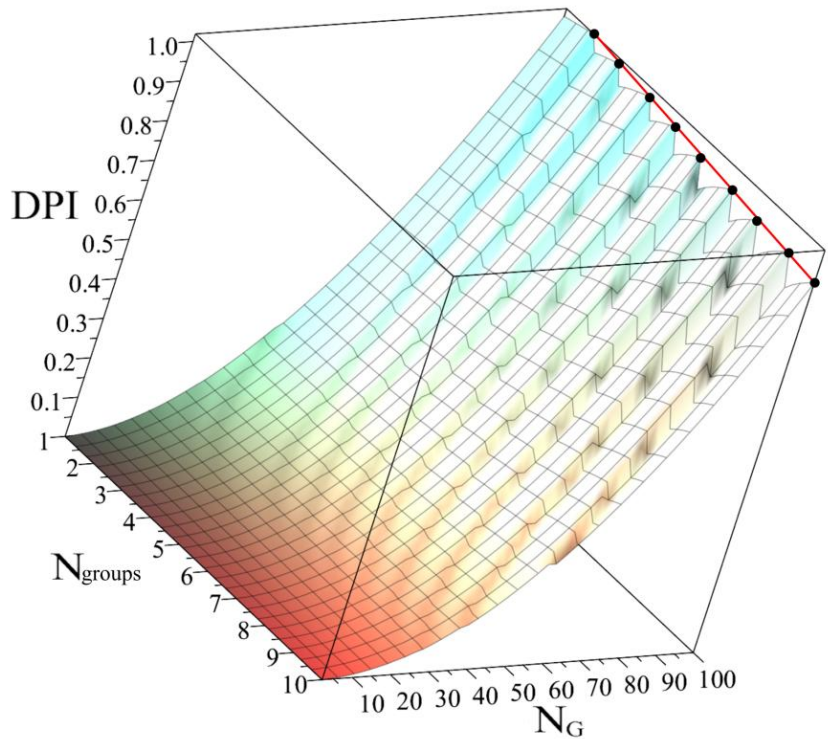

*Fig. 14 Dependence of DPI on $N_{groups}$ - number of groups of synchronized oscillators and $N_G$ – number of synchronized oscillators per group*

The cross section of the chart from Fig.(14) for the four cases of first type group numbers $N_{Groups1}$ =1,4,7,10 first type groups and the second engaging other oscillators not involved in the first structure is presented in Fig(15). Similarly to presented earlier cases, dependence is growing up together with growth of size of the structures. Note that in the case of $N_{Groups1}$ =1, DPI is achieving finally value 1 for the case of 100 synchronized oscillators in first group, meaning complete synchronization of the whole system. This is the only case of complete synchronization as it demands value $N_{Groups1}$ =1.

In the Fig.(15) there is also presented window with focused part of the chart. It presents uniqueness of DPI which is not obvious in this section while looking for the whole range.

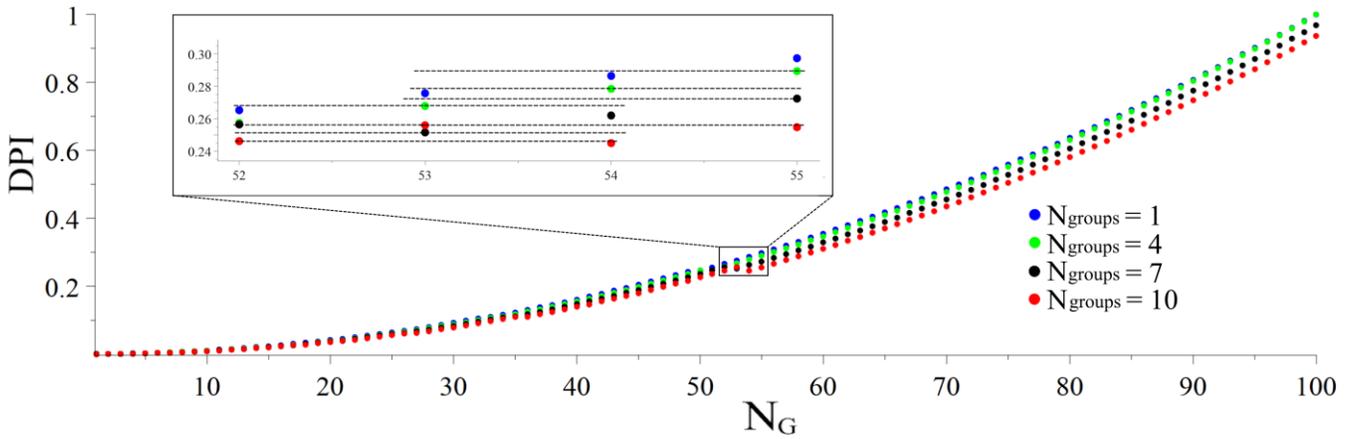

*Fig. 15 Dependence of DPI on number $N_G$ of synchronized oscillators in the first group structure - the cross section for the case of constant, one existing first type group.*

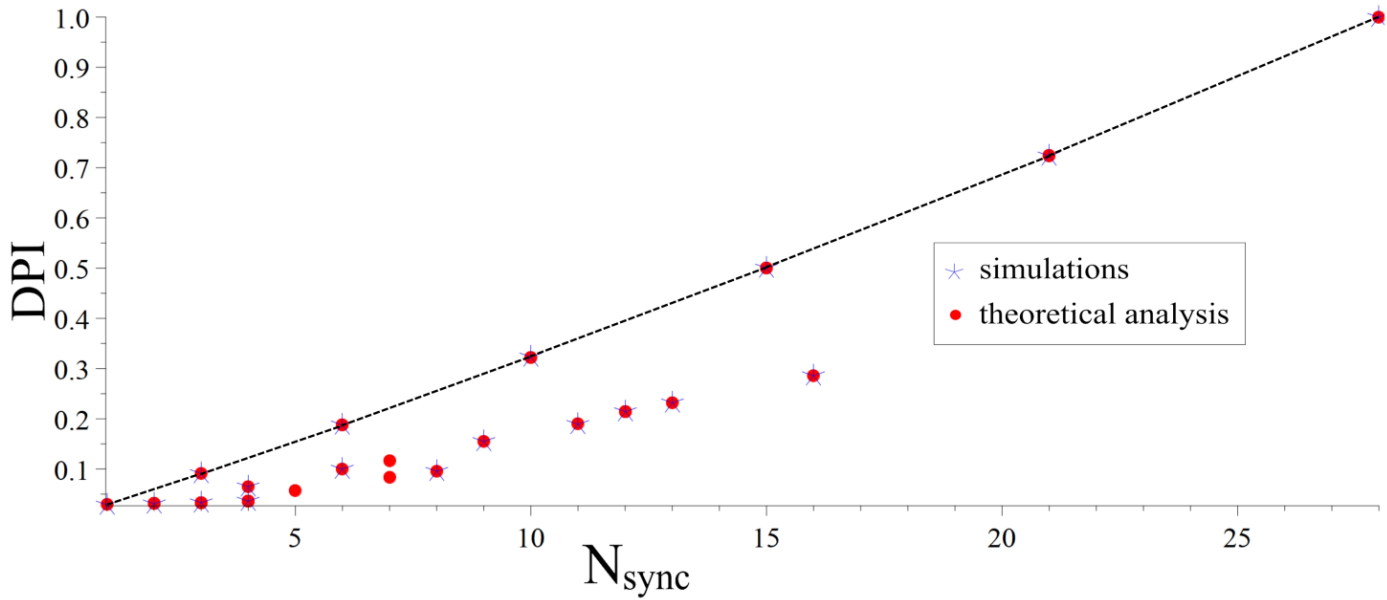

*Fig. 16 Comparison of dependence of DPI on $N_{sync}$ - number of pairs of synchronized oscillators for obtained during investigations and analytical analysis of numbers of synchronized oscillators per group in case of 8 coupled oscillators.*

In the Fig.(16) there is presented comparison of DPI values obtained from our investigations marked with stars and theoretical analyses marked with dots. DPI is shown as dependent on number of synchronized pairs for the case of 8 coupled oscillators. The first important aspect is that generally together with growing number of synchronized pairs DPI is also growing up to final value equal to 1 for 28 synchronized pairs. It refers to maximal number of two-element combinations with no repetitions from the set of 8 elements which in this case means complete synchronization. The other maximal numbers of two element combinations for 2-7 oscillators are connected with black dashed line.

In the Fig.(16) there can be also seen three red dots without respective star points. It means that only 3 configurations from all potentially possible were not detected during our investigations. We suspect that these configurations would appear for other ranges of parameters or initial conditions. One way or another, generally presented results prove that analysed systems do not include in their dynamics any tendency to direction it in some special configurations other way than statistically equal configurations, dependent on the system parameters and initial conditions.

Note also that appearing in the Fig.(16) different DPI values for the same $N_{sync}$ numbers comes from the fact that $N_{sync}$ is not the only parameter identifying dynamical state of the system.

| | n ALL | n TDLE | Aver Rate | DPI | 1.0 | 0.724 | 0.5 | 0.323 | 0.286 | 0.232 | 0.214 | 0.189 | 0.188 | 0.166 | 0.155 | 0.1 | 0.095 | 0.091 | 0.065 | 0.036 | 0.033 | 0.031 | 0.029 |
|---|---|---|---|---|---|---|---|---|---|---|---|---|---|---|---|---|---|---|---|---|---|---|---|
| 8 OSC | 86 | 63 | 0,73256 | config | 2 8 0 | 1 2 1 1 | 1 1 5 2 | 1 1 0 3 | 1 1 6 0 | 2 1 3 0 | 2 1 2 0 | 2 1 1 1 | 2 6 4 1 | 7 2 2 | 9 1 2 | 6 2 2 | 8 0 2 | 3 5 1 | 4 3 2 | 4 0 4 | 6 2 3 | 2 4 2 | 1 6 1 |
| | | | | n TDLE | 1 | 1 | 4 | 2 | 4 | 5 | 5 | 2 | 5 | 2 | 4 | 3 | 1 | 3 | 7 | 2 | 4 | 4 | 4 |
| | | | | n ALL | 1 | 1 | 4 | 5 | 4 | 5 | 5 | 6 | 6 | 5 | 6 | 6 | 2 | 5 | 2 | 3 | 6 | 10 | 4 |

*Table 1. Results obtained from numerical simulations for 8 coupled oscillators.*

DPI values shown in Fig.(16) are also presented in Table (1) together with the respective configurations ($N_{sync}, N_{unsync}, N_{Groups}$). Beneath these values there are given numbers of these configurations detected by TDLE method and ALL existing.

Note firstly that opposite to results from Fig.(16) all the presented in Table (1) come from numerical investigations which is the cause of not appearing and mentioned above 3 points with no respective numerical solutions.

Note also that appearing in Table (1) different numbers of the same configurations mean identified symmetrically equivalent states. It can be seen from Table (1) that not all of these symmetrical states were detected by proposed TDLE method. As analysed beneath the reason was the way of obtaining introduced threshold function which was identified for the simplest case of four coupled oscillators. Thus for more complex case of 8 coupled oscillators the Average Rate of properly detected states for considered case is equal to 0.73 not 1. Beneath we propose easy solution for future improvement of our algorithm.

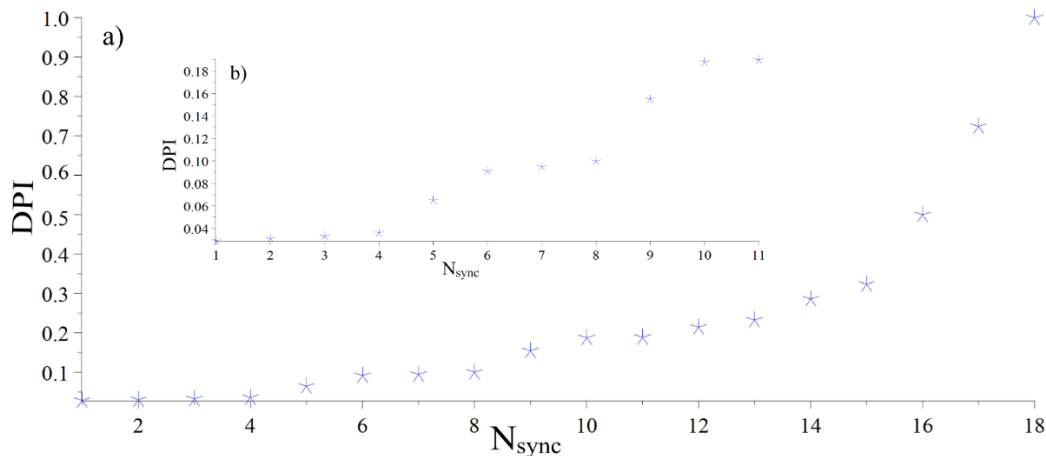

*Fig. 17 Dependence of DPI on growing up synchronized structures obtained during investigations for the case of 8 coupled oscillators. Full range (23 a); uniqueness of DPI for focused range (23.b).*

As there is no existing a single parameter depicting growth of complexity of dynamics of coupled oscillators, what was mentioned above while analysing Fig (16), to present monotonical growth of DPI together with growth of synchronized structures, presented results obtained during simulations for the case of 8 coupled oscillators were rearranged. They were ordered depending on DPI values and results can be seen in Fig.(17). Fig.(17a) is showing general growth of DPI values and its uniqueness. For better presentation of the uniqueness in the range of lower DPI values they are presented in focused view in Fig.(17b).

## Table 2

| | n ALL | n TDLE | Aver Rate | DPI | 1.0 | 0.724 | 0.5 | 0.323 | 0.286 | 0.232 | 0.214 | 0.189 | 0.188 | 0.166 | 0.155 | 0.1 | 0.095 | 0.091 | 0.065 | 0.036 | 0.033 | 0.031 | 0.029 |
|---|---|---|---|---|---|---|---|---|---|---|---|---|---|---|---|---|---|---|---|---|---|---|---|
| **8 OSC** | 86 | 63 | 0,73256 | config | 28;0;1 | 21;1;1 | 15;2;1 | 10;3;1 | 16;0;2 | 13;0;2 | 12;0;2 | 11;1;2 | 6;4;1 | 7;2;2 | 9;1;2 | 6;2;2 | 8;0;2 | 3;5;1 | 4;3;2 | 4;0;4 | 6;2;3 | 2;4;2 | 1;6;1 |
| | | | | n TDLE | 1 | 1 | 4 | 2 | 4 | 5 | 5 | 2 | 5 | 2 | 4 | 3 | 1 | 3 | 7 | 2 | 4 | 4 | 4 |
| | | | | n ALL | 1 | 1 | 4 | 5 | 4 | 5 | 5 | 6 | 6 | 5 | 6 | 6 | 2 | 5 | 2 | 3 | 6 | 10 | 4 |

| | n ALL | n TDLE | Aver Rate | DPI | 0,68 | 0,44 | 0,26 | 0,25 | 0,21 | 0,16 | 0,14 | 0,12 | 0,09 | 0,08 | 0,08 | 0,04 |
|---|---|---|---|---|---|---|---|---|---|---|---|---|---|---|---|---|
| **7 OSC** | 50 | 30 | 0,6 | config | 15;1;1 | 10;2;1 | 11;0;2 | 6;3;1 | 9;0;2 | 8;1;2 | 6;1;2 | 3;4;1 | 4;2;2 | 4;3;2 | 5;0;3 | 1;5;1 |
| | | | | n TDLE | 1 | 3 | 3 | 3 | 3 | 6 | 2 | 4 | 0 | 2 | 0 | 3 |
| | | | | n ALL | 1 | 3 | 3 | 5 | 3 | 7 | 5 | 4 | 9 | 6 | 1 | 3 |

| | n ALL | n TDLE | Aver Rate | DPI | 0,23 | 0,2 | 0,17 | 0,13 | 0,07 | 0,06 | 0,05 |
|---|---|---|---|---|---|---|---|---|---|---|---|
| **6 OSC** | 26 | 24 | 0,92308 | config | 7;0;2 | 6;0;2 | 3;2;1 | 4;1;2 | 3;0;3 | 2;2;2 | 1;4;1 |
| | | | | n TDLE | 3 | 3 | 3 | 4 | 3 | 5 | 3 |
| | | | | n ALL | 3 | 3 | 3 | 5 | 3 | 6 | 3 |

| | n ALL | n TDLE | Aver Rate | DPI | 0,55 | 0,25 | 0,2 | 0,09 | 0,08 |
|---|---|---|---|---|---|---|---|---|---|
| **5 OSC** | 10 | 10 | 1 | config | 10;1;1 | 3;2;1 | 4;0;2 | 2;1;2 | 1;3;1 |
| | | | | n TDLE | 1 | 2 | 2 | 3 | 2 |
| | | | | n ALL | 1 | 2 | 2 | 3 | 2 |

| | n ALL | n TDLE | Aver Rate | DPI | 0,43 | 0,17 | 0,13 |
|---|---|---|---|---|---|---|---|
| **4 OSC** | 5 | 5 | 1 | config | 3;1;1 | 2;0;2 | 1;2;1 |
| | | | | n TDLE | 1 | 2 | 2 |
| | | | | n ALL | 1 | 2 | 2 |

*Table 2. Results obtained from simulations for 4-8 coupled oscillators.*

Results of DPI values obtained from simulations for all the investigated CN consisted of 4…8 coupled oscillators are presented in Table (2). Similarly to Table (1) in the line below DPI, the respective configurations ($N_{sync}, N_{unsync}, N_{Groups}$) are presented. Beneath configurations there are given numbers of configurations detected by TDLE method and ALL existing. It is also presented Average Rate of properly detected states for considered cases. Note that results are ordered

Within these results it can be seen that for the cases of 4 and 5 coupled oscillators all the configurations were properly detected and Average Rate is equal to one. Simultaneously Average Rate lowers down together with growth of number of coupled oscillators. It proves that crucial aspect in improvement of efficiency of the presented approach is finding simple universal procedure of threshold determination. As mentioned earlier to simplify investigations the threshold was obtained for the smallest CN consisted of 4 coupled oscillators. In the next stage of our investigations we are planning to find an efficient way of threshold normalization that it would be more universal and appropriate for more complex cases.

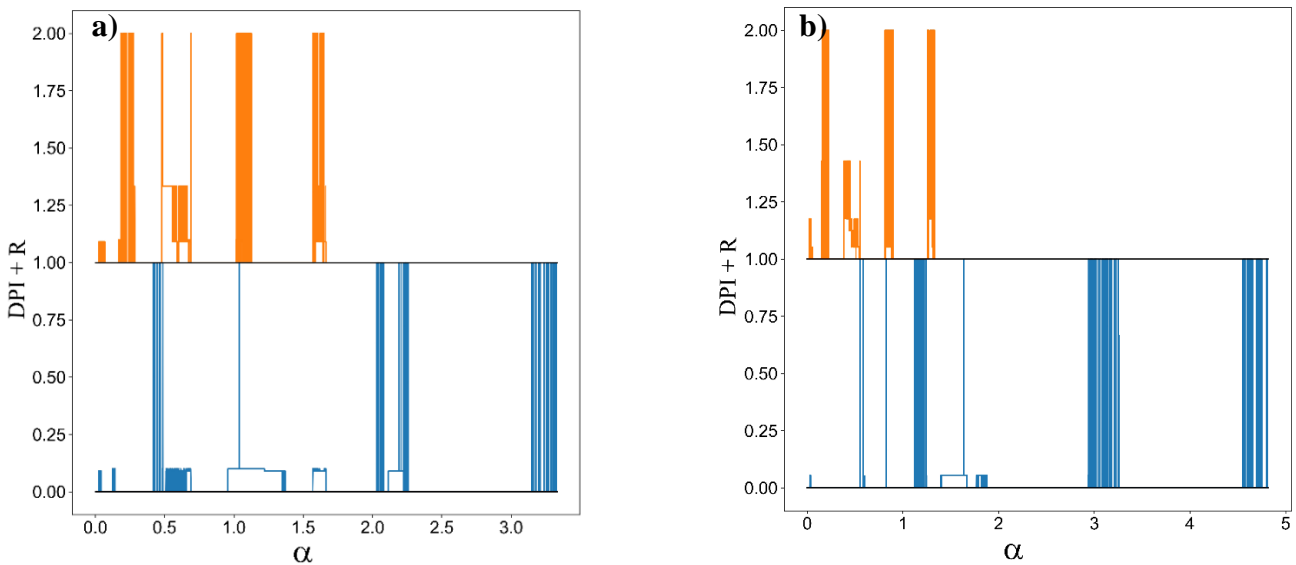

*Fig. 18 Dependence of DPI on coupling coefficient α and radius of coupling R={1,2} for 4 oscillators (20 a) and 5 oscillators (20 b)*

Another feature of the proposed DPI index comes from its normalized values achieving values in the range between 0 and 1. It allows to compare dynamical phenomena for different sizes of the networks of coupled oscillators by adding these values to DPI. Te same way it can be also extended on other information for instance radius of coupling and this way preventing easy localization of system topology for which dynamical phenomena was detected. Examples are presented in Figs(18 - 20). Comparing these figures different aspects of influence of system complexity, radius of coupling, and coupling coefficient on system's dynamics can be taken.

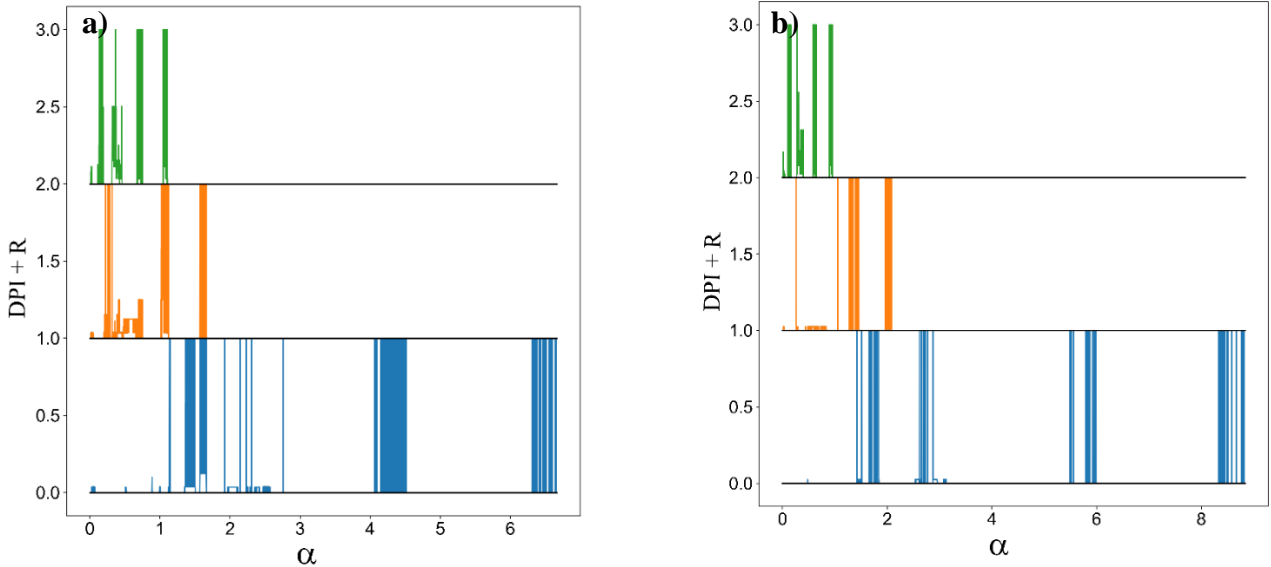

*Fig. 19 Dependence of DPI on coupling coefficient $\alpha$ and radius of coupling R={1,2,3} for 6 oscillators (21 a) and 7 oscillators (21 b)*

Firstly it can be seen that together with growth of the radius of coupling the systems show a greater tendency to bifurcate between different dynamical states while changing coupling coefficient. These increasing numbers of bifurcations between different final configurations are visible in all the figures proving independence of this behaviour on number of coupled oscillators.

The second observation is that more complex states meaning bigger values of DPI index different from 1 (the case of complete synchronization), are also connected with radius of coupling and exceed value 0.5 only for radius equal to 4 in the case system of 8 coupled oscillators presented in Fig.(20).

This fact together with earlier conclusion taken while analysing the Fig.(16) and suggesting existence of only statistically equivalent solutions, leads to following conclusion: for the case of more complex networks of Duffing oscillators the most probable is existence of more complex solutions but not having special regularities worth to investigate.

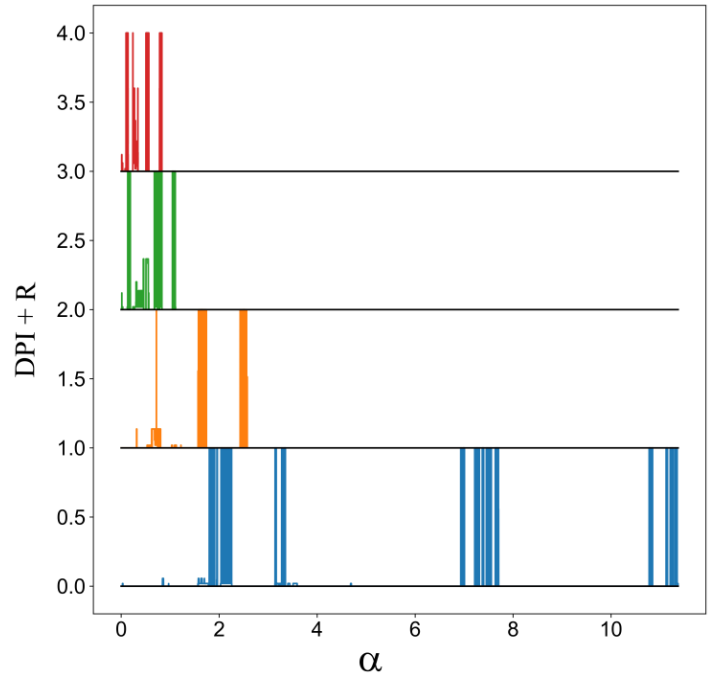

*Fig. 20 Dependence of DPI on coupling coefficient $\alpha$ and radius of coupling R={1,2,3,4} for 8 oscillators*

**Results of time efficiency of presented approach**

Time efficiency of the proposed algorithm is presented in Fig.(21). The rate given in this figure was calculated as CDP detection time of the presented method compared to the CDP detection time of the perturbation vector analysis algorithm (time at which the first pair of oscillators synchronize regarding perturbation decay with changes within $10^{-5}$).

Results of time efficiency comparison for all investigated topologies of CNS are also summarized in Table (3). As we can see, for each investigated threshold function, the average time rate is in the range of 30%, indicating that the proposed algorithm is approximately three times faster than the perturbation vector analysis algorithm. The worst result was obtained for R73 topology, where the time rate varied from 40% to 50% across different threshold functions. This indicates that, even in the worst-case scenario, the proposed algorithm detects complex dynamical phenomena approximately twice as fast as the perturbation vector analysis algorithm. Note, that different threshold functions slightly affect the detection time of the proposed algorithm.

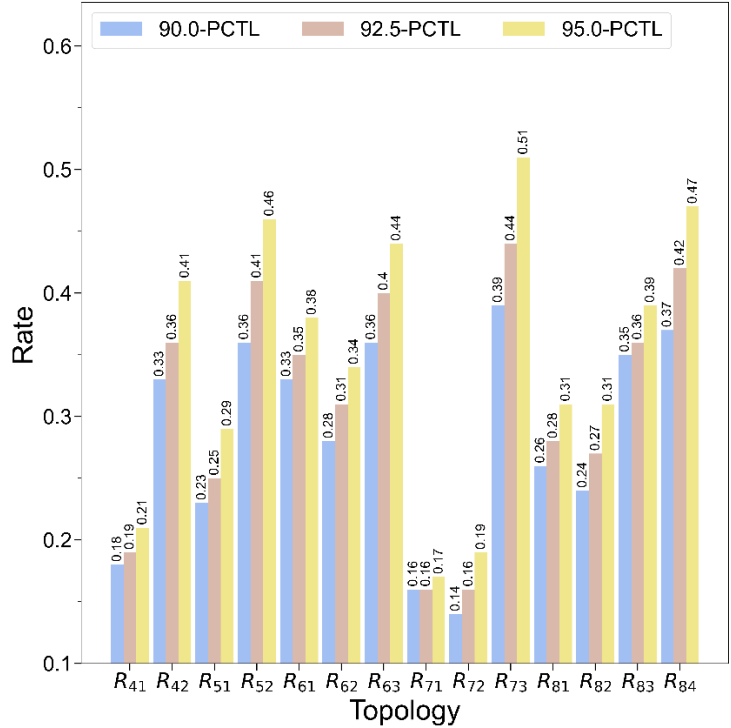

*Fig. 21 Time efficiency of the algorithm*

| Threshold function | Average efficiency-case[%] | Lowest efficiency-case[%] |
|---|---|---|
| 90- PCTL | 28.46 | 39.21 |
| 92.5- PCTL | 31.13 | 44.35 |
| 95- PCTL | 34.88 | 50.78 |

*Table 3 Time efficiency*

## 4. Conclusions

In the article we proposed new approach for fast scanning of CN parameters and initial conditions to obtain a preliminary approximation of its dynamic state.

We utilized presented in our earlier article TDLE index and its feature of early signalizing of spreading CN dynamics into synchronized subsystems. We demonstrated that this approach is highly effective tool while analysing complex phenomena existing in networks of coupled nonlinear systems. Based on this approach we also introduced new DPI index depicting dynamical state of CN. This index was deeply analytically studied and compared with results of numerical simulations confirming correctness of its assumed features. While we proved DPI uniqueness it can be used not only to depict dynamical state of CN but also to recover its dynamical configuration in cases where some interesting phenomena exist.

Based on the presented results we took general conclusions on dynamical state of analysed CN.

The most important fact is that proposed approach is highly universal and can be applied for both, symmetrical and non-symmetrical topologies of coupling as well as networks of identical and non-identical oscillators.

Moreover, since DPI values are obtained from the current state of dynamical system given by values of the system variables, proposed method of fast searching has a huge potential for experimental application.
.

………………………………………………………………

**ACKNOWLEDGMENTS**

The author declares that has no conflict of interests.

All the data will be made available on reasonable request.

V.D has been supported by the National Science Centre, Poland, PRELUDIUM Program (Project No. 2023/49/N/ST8/02436) This paper has been completed while the first author was the Doctoral Candidate in the Interdisciplinary Doctoral School at the Lodz University of Technology, Poland.